\documentclass[11p,reqno]{amsart}
\usepackage{amsmath}
\usepackage{amsfonts}
\usepackage{amssymb}
\usepackage{graphicx}
\usepackage{color}
\newtheorem{theorem}{Theorem}[section]
\newtheorem*{theorem*}{Theorem}
\newtheorem{lemma}{Lemma}[section]
\newtheorem{corollary}[theorem]{Corollary}
\newtheorem{proposition}{Proposition}[section]

\newtheorem{remark}[theorem]{Remark}

\def\Ric{\text{Ric}}

\def\Ric{\operatorname{Ric}}
\def\Rm{\operatorname{Rm}}

\def\tr{\operatorname{tr}}

\def\Rm{{\operatorname R}}
\def\Ric{{\operatorname{Ric}}}

\def\ad{{\operatorname{ad}}}

\numberwithin{equation}{section}

\begin{document}
	\title[Holonomy of complex manifolds]{Holonomy and the Ricci curvature of complex Hermitian  manifolds}

\author{Lei Ni}
\address{Zhejiang Normal University, Jinghua, Zhejiang, China, 321004;
Department of Mathematics, University of California, San Diego, La Jolla, CA 92093, USA}
\email{leni@zjnu.edu.cn; lni.math.ucsd@gmail.com}


\subjclass[2010]{}

\begin{abstract} We prove two results on geometric consequences of the representation of  restricted holonomy group of a Hermitian connection. The first result concerns when such a Hermitian manifold is K\"ahler in terms of the torsion and the irreducibility of the holonomy action. As a consequence we obtain a criterion of when a Hermitian manifold (and connection) is a generalized Calabi-Yau (in the sense that the Chern Ricci vanishes or equivalently that the restricted holonomy is inside $\mathsf{SU}(m)$). The second result concerns when a compact K\"ahler manifold with a generic restricted holonomy group is projective.
\end{abstract}

\maketitle

\section{Introduction}

The holonomy group of the Levi-Civita connection (denoted by $D$), which is called the Riemannian holonomy,  is an important object for the study of Riemannian manifolds \cite{Berger, Besse, Kostant-TAMS}. Roughly the size of the group measures how much the manifold locally  is deviated from being Euclidean.

 For complex manifolds equipped with a Hermitian metric. There are several canonical Hermitian connections, including the Chern connection (denoted as $\nabla^c$ satisfying $(\nabla^c)''=\bar{\partial}$), the Lichnerowiz connection ($=\frac12 (D-JDJ)$), the Bismut/Strominger connection (denoted by $\nabla^b$, which is characterized as its torsion being a $3$-form), and a family of Gauduchon connections ($\nabla^t=(1-\frac{t}{2})\nabla^c+\frac{t}{2}\nabla^b$, which consists of the linear combination of the Chern connection and the Bismut connection). Here the connection $\nabla$ is called Hermitian if it keeps the metric and the associated almost complex structure $J$ invariant, namely satisfies $\nabla g=0, \nabla J=0$.  The general study of the interaction between the holonomy group of Hermitian connections and the associated geometries of these Hermitian connections is scarce in the literature. Recently in \cite{AndV} the holonomy of Bismut connection on Hopf manifolds endowed with the Bismut connection was computed.

 The holonomy group of these canonical connections on a complex Hermitian manifold seems quite different from the Riemannian holonomy.  The study of these connections is complicated by the lack of a de Rham decomposition theorem which splits the Riemannian manifold via the Riemannian holonomy actions. From the  recent papers \cite{NZ-TAMS, NZ-CAS2}, one can infer that for a class of Hermitian homogenous manifolds, namely the locally Chern homogenous Hermitian manifolds, the restricted holonomy groups are the same as the ones of Hermitian symmetric spaces  of  lower dimensions if the Hermitian connection is not K\"ahler. In fact this follows from a classification result asserting that {\it A  simply-connected locally Chern homogenous Hermitioan manifold must be product of a complex Lie groups (equipped with a left invariant Hermitian metric and the Chern connection which makes a set of  left invariant frame unitary and invariant under the parallel transport \cite{Boothby}), and a Hermitian symmetric space. Hence if the manifold is irreducible, it is either a complex Lie group or a Hermitian symmetric space.} The Chern connection of the above metric  on the complex Lie group has a trivial holonomy $\{\operatorname{id}\}$.

The following results of Cartan and Ambrose-Singer \cite{AS, Kostant-Nagoya} (as well as its Hermitian analogue) provide the background and motivations for the above mentioned study  in \cite{NZ-TAMS, NZ-CAS2}. Even though the statements of the theorem blow  do not involve the holonomy groups of the corresponding connections, the proof however used it \cite{Helgason, Kostant-Nagoya}.

\begin{theorem}[Cartan, Ambrose-Singer]\label{Cartan-etc}
Let $(M^n, g)$ be a simply-connected complete Riemannian manifold. Then the statements (a) and (b) are equivalent for cases (i) and (ii).

(i) (a) $(M^n, g)$ is a Riemannian symmetric space; (b) The Riemannian curvature tensor  $\Rm^{D}$ is parallel with respect to $D$, i.e. $D \Rm^D=0$.

(ii) (a) $(M^n, g)$ is Riemannian homogeneous (in the sense that the isometry group of $M$ acts transitively on $M$); (b) There exists a metric connection $\nabla$ such that its curvature $\Rm^\nabla$ and torsion tensor $T^\nabla$ are parallel with respect to $\nabla$.
\end{theorem}

To prove the main results  in \cite{NZ-TAMS, NZ-CAS2}, studying  the holonomy (of the Chern connection) action on $T^{1, 0}M$ plays an important role. In fact as a step towards the main classification theorem, it was proved that {\it for a complex Hermitian manifold $(M^m, g)$ (with $m=\dim_\mathbb{C}(M)$) whose Chern connection is Ambrose-Singer, if the universal cover does not contain any  de Rham factor which is K\"ahler, then the Chern Ricci curvature of $M$ must be flat. In particular its restricted holonomy must be inside $\mathsf{SU}(m)$.}
This result partially answers a question of Yau \cite{Yau-op} where he asked if {\it  one can say something nontrivial about the compact Hermitian manifolds whose holonomy is a proper subgroup of $\mathsf{U}(m)$}, by asserting that for a Chern homogeneous manifold, the holonomy of the Chern connection reduces to the holonomy of its Hermitian symmetric space factor.  Recall that a Hermitian connection $\nabla$ is called Ambrose-Singer if both its curvature and torsion are parallel with respect to $\nabla$. By the theorem of Ambrose-Singer/Kostant \cite{Kostant-Nagoya}, these conditions imply that the manifold $(M^n, g)$ with such a Hermitian connection admits a transitive holomorphic-isometrical action of a connected Lie group $G$,  which also keeps the connection $\nabla$ invariant.

Below is  a result which also provides a restriction of the holonomy group of  a{\it (any)} Hermitian connection $\nabla$ of a Hermitian manifolds $(M^m, g)$.

\begin{theorem} \label{thm:CY-torsion} For a complex Hermitian manifold $(M^m, g)$, and a Hermitian connection $\nabla$, assume that  the restricted holonomy group $G$ action on $T^{1, 0}(M)$ is irreducible. Assume  further that $\Ric^{1}\ne  0$ and the torsion of $\nabla$ is parallel. Then $(M^m, g)$ must be K\"ahler.  Equivalently, for any non-K\"ahler Hermitian manifold whose restricted holonomy with respect to a Hermitian connection with a  parallel torsion is irreducible,  it must be a generalized Calabi-Yau in the sense that its  Chern Ricci vanishes.
\end{theorem}

Recall that $\Ric^{1}(X, Y):=\sum_{i=1}^m \langle \Rm_{X, Y}(e_i), \overline{e_i}\rangle$, with $\{e_i\}_{i=1}^m$ being a unitary frame of $T^{1, 0}_pM$, is the Chern Ricci curvature (also called the first Ricci curvature).
The result is a bit surprising given  that it holds  for any  Hermitian connections, while the results of \cite{NZ-TAMS, NZ-CAS2} hold only for the Chern connection. Given that for K\"ahler manifolds, the torsion, being zero, is parallel, and the generic K\"ahler metrics have irreducible holonomy, the assumptions of the theorem do not seem too unnatural.  The above result puts a restriction on the restricted holonomy of a Hermitian connection with a parallel torsion on a non-K\"ahler manifold by asserting that it is either reducible or a subgroup of  $\mathsf{SU}(m)$,   answering partially again the above-mentioned Yau's question for all Hermitian connections with parallel torsion, without assuming  the compactness of the manifold. Under the balanced condition of the metric, the Chern-Ricci flatness was proved earlier for Gauduchon connections in \cite{Vez}.  Given that the study of non-K\"ahler manifolds (and non-K\"ahler Calabi-Yau) with special properties of its torsion has been found necessary in string theory (cf. \cite{GHR, Stro}), the above result may have some physical implications as well.

There are  recent studies/results which assert,  for  a Hermitian manifold under some assumptions,  either for the Bismut connection to conclude  its curvature tensor must be K\"ahler-like (namely satisfies all the symmetries of the curvature of a K\"ahler metric) (cf. \cite{AOUV, Zhao-Zheng}), or for the Gauduchon connection $\nabla^t$ to conclude that it must be  K\"ahler \cite{LS}.
More precisely, in \cite{Zhao-Zheng} it was proved that for a Bismut connection  {\it  the  associated curvature tensor is K\"ahler-like if and only if the torsion is parallel and the metric is pluriclosed (namely $\partial\bar{\partial} \omega=0$ with $\omega$ being the K\"ahler form).} In \cite{LS} it was proved that  for a Hermitian manifold $(M^m, g)$,  {\it if the curvature tensor of  the Gauduchon connection $\nabla^t$  is K\"ahler-like (or in particular, flat),  and the norm square of the Gauduchon one-form can attain its maximum somewhere,  then $g$  is K\"ahler if $t\ne 0, \frac23, \frac45, 2$}. In addition to these result, Theorem 1.10 of a more recent paper \cite{NZ-CAS2}  asserts that {\it for  a compact Hermitian manifold whose Gauduchon connection $\nabla^t$ (with $t\ne 0$) is Ambrose-Singer, and with vanishing first and third Ricci curvatures, it  must be K\"ahler.}

In view of that Theorem \ref{thm:CY-torsion} appears differently from the above mentioned three results, it is not surprising that the proof here is also orthogonal to the proofs of above mentioned results, which  involves careful analysing the derivatives of torsion and curvature tensors. The proof here appeals to some simple considerations via the Lie groups/algebras.

The second result of the paper concerns the projectivity of a compact K\"ahler manifold, namely when a compact K\"ahler manifold can be embedded into some complex projective space of a higher dimension. Since we are on a K\"ahler manifold, the holonomy group involved is the well-studied Riemannian holonomy \cite{Besse}, and there is only one Ricci curvature.  The projectivity (which is the same as being an algebraic manifold) of a compact complex manifold was characterized by the existence of a K\"ahler form whose cohomology class is integral \cite{MK}, which is equivalent to the existence of a positive holomorphic line bundle. There are many K\"ahler manifolds which are not projective. Some could not even be deformed into a projective one since that would imply that they have the homotopy type of projective manifold. It is interesting to characterize projectivity using intrinsic properties of a K\"ahler manifold.
There are  recent works (\cite{N}, \cite{NZ}, \cite{NZ-2}, \cite{Yang})  asserting  projectivity of a compact K\"ahler manifold $(M, g)$ under  various positivity assumptions including the positivity of the holomorphic sectional curvature $H$ (cf. \cite{Yang}), the positivity of $\Ric^\perp$ (cf. \cite{NZ}),  the positivity of the second scalar curvature $S_2$ (cf. \cite{NZ-2}), and the  positivity of $\Ric_k$ (or even $2$-positivity of $\Ric_k$, $\Ric_3^\perp$), $\Ric^+$, $S^\perp_2$ \cite{N}, or even the much weaker $\mathcal{BC}$-2 positivity (cf. Section 4 of \cite{N}). These results and  positivity are also related to the rational connectedness of $M$ \cite{CDP, N, Yang}. We refer  the interested readers to the above papers for the detailed definitions of these curvature conditions and the statements of the results.

On the other hand, projectivity was also proved for compact Calabi-Yau manifolds (namely with a K\"ahler metric of $\Ric\equiv 0$ and restricted holonomy being $\mathsf{SU}(m)$) by appealing to the Lie/representation theory \cite{Beau} (cf. \cite{Besse, Wu} for expositions). There the holonomy group \cite{Berger} and its representation hold the key. Since $\Ric\equiv 0$ is a special case of $\Ric\ge 0$, one natural additional condition needed is the nonnegativity of some curvatures, including $ H, \Ric_k, \Ric^{+}, \Ric^\perp, S_2$ discussed in the above paragraph. The second goal of this paper is the following  result concerning the relationship between the holonomy group, the projectivity and nonnegativity of the curvature (which includes the $\Ric\equiv 0$ case).

\begin{theorem}\label{thm:main1}Let $(M^m, g)$ be a compact K\"ahler manifold with the restricted holonomy being  $\mathsf{U}(m)$ (or $\mathsf{SU}(m)$). Assume   any one of the following:

(i) The sum of the smallest two eigenvalues of  $\Ric$ is nonnegative;

(ii) The mixed curvature $\mathcal{C}_{\alpha, \beta}(X)\ge 0$, for $\alpha> 0, 3\alpha+2\beta> 0$;\footnote{When we were writing up this result, we were informed that in a recent ArXiv preprint \cite{CL}, a  result  close to part (ii) below appeared. The argument of \cite{CL} uses a conformal change, while the proof here follows closely the argument of \cite{Tang}. In \cite{Tang} the projectivity in \cite{NZ-2} under the positivity of $S_2$ was extended to the quasi-positive case. }

(iii) The second scalar curvature $S_2\ge 0$.

Then $M$ is projective and rationally connected (for $\mathsf{U(m)}$ holonomy).
\end{theorem}

The mixed curvature $\mathcal{C}_{\alpha, \beta}(X)$ was introduced in \cite{CLT} to extend results under the positivity assumption of $\Ric_k$, for some $1\le k\le m$ (with $k=1$ being the holomorphic sectional curvature),   $\Ric^\perp$ or $\Ric^+$ in \cite{Yang, NZ} and \cite{N}. By the definition $\mathcal{C}_{\alpha, \beta}(X)\ge 0$ if and only if
\begin{equation}\label{eq:mixed}
\mathcal{C}_{\alpha, \beta}(X)=\alpha|X|^2 \Ric(X, \overline{X})+\beta \Rm(X, \overline{X}, X, \overline{X})\ge 0, \forall X \in T^{1,0}_xM, \forall x\in M.
\end{equation}
For example $\alpha=1$, $\beta=-1$ corresponds to $\Ric^\perp$.

The results of \cite{Yang, NZ, N, CLT} were extended to the quasi-positive case recently for holomorphic sectional curvature in \cite{Zhang-Zhang} and mixed curvature, as well as $S_2$, in \cite{Tang}. As pointed out in \cite{CLT} (cf, Lemma 2.1 therein), part (ii) includes the case $\Ric_k\ge0$, $\Ric^{\perp}\ge 0$ and $\Ric^+\ge 0$. Hence one can view the above result as the rigidity cases of the previous results of \cite{NZ, N, CLT, Zhang-Zhang, Tang}.   In \cite{N}, the projectivity was proved for K\"ahler manifold with the $\mathcal{BC}$-2 positivity,  namely if $\forall x\in M$, for any two dimensional subspace $\Sigma\subset T^{1,0}_xM$ there exists a vector  $v\in T_x^{1, 0}M$ such that
$$
\int_{\mathbb{S}^3} \Rm(v, \bar{v}, Z, \bar{Z})\, d\theta(Z)>0.
$$
This is formally  weaker than $S_2>0$ since one can easily derive the above if $S_2>0$. It is interesting to find out whether or not the projectivity result under  the $\mathcal{BC}$-2 positivity in \cite{N} still holds when  the above `$>0$' is replaced by `$\ge0$' together with a generic holonomy assumption (as well as if the positivity can eb replaced by the quasi-positivity). This question applies similarly to other positive curvature conditions considered in \cite{N}, but not being covered by the above theorem, namely the $S_2^\perp\ge 0$, $\Ric_k$ 2-nonnegative and $\Ric^\perp_3\ge 0$ cases. \footnote{Very recently  Y. Niu and her collaborator have successfully pushed the arguments to cover these three cases.}

As in the study of the Ricci flow invariant cone conditions \cite{Wilking}, the  rigidity for the curvature conditions in the boundary of the cones  is expected to be related to the holonomy group of the manifold. The above theorem partially supports this expectation from a different perspective.

The last result is an observation that for any Hermitian manifold whose Bismut connection is with the parallel torsion, $T_{\mathbb{C}}M$ is a complex Lie group bundle with the Lie bracket being defined by the torsion of the Chern connection. Moreover the splitting $T_{\mathbb{C}}M=T^{1, 0}M\oplus T^{0,1}M$ is also a direct sum of ideals.

\section{Preliminaries}
Let $G$ be the restricted holonomy of a metric connection. Since the restricted holonomy group with respect to a metric connection is a path-connected subgroup of $\mathsf{SO}(n)$, hence it is a Lie subgroup. Moreover it is closed if it acts irreducibly by Theorem 2 on page 277 of \cite{KN}. Therefore it is a compact Lie group. We denote its Lie algebra by $\mathfrak{g}$. Let $B$ be the Cartan-Killing form of $\mathfrak{g}$. It is well known that  there exists a $\ad_{\mathfrak{g}}$-invariant Euclidean inner product on $\mathfrak{g}$. In our case since $G$ is a compact subgroup of $\mathsf{SO}(n)$, one can easily check that for any $x, y\in \mathfrak{so}(n)$
$$
\langle x, y\rangle =-\operatorname{tr}(x y)=\sum_{i=1}^n \langle x(\epsilon_i), y(\epsilon_i)\rangle=\operatorname{tr}(x {}^ty)
$$
is an $\ad$-invariant positive definite metric on its Lie algebra $\mathfrak{so}(n)$, where $\{\epsilon_i\}_{i=1}^n$ is  an orthonormal frame of $V=T_pM$ (when $M$ is a Riemannian manifold). Here ${}^ty$ is the transpose of the transformation (matrix) $y$. This {\it extrinsic} invariant inner product is important for our proof since we need to consider the inner product of $\mathfrak{g}$ with an element of $\mathfrak{so}(n)$ which possibly is not inside $\mathfrak{g}$.

For a Hermitian manifold, we can choose a frame consists of $\{\epsilon_i, J\epsilon\}_{i=1}^m$ as an orthonormal frame adapted to the almost complex structure $J$.  Since the parallel transport preserves $J$, $G$ is a compact subgroup of $\mathsf{U}(m)\subset \mathsf{SO}(2m)$. The inner product naturally extends to its Lie algebra $\mathfrak{u}(m)$ as
$$
\langle a, b\rangle =-2\Re \operatorname{tr}(a b)=2\Re\sum_{i=1}^m \langle a(e_i), \overline{b(e_i)}\rangle
$$
where $\{e_i\}_{i=1}^m$,  with  $e_i=\frac{1}{\sqrt{2}}(\epsilon_i-\sqrt{-1}J\epsilon_i)$, is an unitary frame of $T^{1, 0}_pM$ when $(M^m, g)$ is a complex manifold equipped with a Hermitian metric $g$, namely $g$ (or the inner product) satisfies  $\langle Jx, Jy\rangle =g(Jx, Jy)=g(x, y)=\langle x, y\rangle$ with $J$ being the naturally induced almost complex structure on $T_pM$.
If we express $a=x+\sqrt{-1}y$, $b=x'+\sqrt{-1}y'\in \mathfrak{u}(m)$ in term of matrices of reals, then $x+{}^tx=0=x'+{}^t x'=0$, $y={}^ty, y'={}^ty'$. Here ${}^tx$ denotes the transpose of $x$.
Direct calculation shows that
$$
\langle a, b\rangle =2\operatorname{tr}(x {}^t x'+y{}^ty').
$$
It is easy to check that the inner product is also $\ad$-invariant.

We  need the following result, which is essentially the Schur's Lemma.

\begin{proposition}\label{prop:schur} Let $G$ be a subgroup of $\mathsf{O}(n)$ which acts irreducibly on $\mathbb{R}^n$. Let $a, b$ be linear transformations of $\mathbb{R}^n$ commute with $G$. Then

(i) If $a$ is symmetric, then $a=\lambda \operatorname{id}$ for some $\lambda \in \mathbb{R}$;

(ii) If $a$ is skew-symmetric then $a=\lambda J$  some $\lambda \in \mathbb{R}$, with $J$ being an almost complex structure, namely $J^2=-\operatorname{id}$;

(iii) If $a\ne 0, b$ are skew-symmetric, and $ab=ba$, then $b=\lambda a$.
\end{proposition}
\begin{proof}
See page 278 of \cite{KN}, Lemma 2 and its proof.
\end{proof}

Another well-known result on compact Lie algebras needed for the proof is the following
\begin{proposition}\label{prop:com-lie-alg}
 The Lie algebra $\mathfrak{g}=\mathfrak{z}\oplus \mathfrak{g}_1\oplus \cdots\oplus \mathfrak{g}_s$ as direct (orthogonal with respect to any adjoint invariant inner product such as the one explicitly given above)  sum of simple ideals, with $\mathfrak{z}$ being the center (namely the maximal ideal $\mathfrak{i}$ satisfying $[x, \mathfrak{g}]=0, \forall x\in \mathfrak{i}$). Any ideal of $\mathfrak{k}$ can be written as a partial sum of the summands. In particular, the derived algebra $\mathfrak{g}'=[\mathfrak{g}, \mathfrak{g}]$ is perpendicular to the center $\mathfrak{z}$.
\end{proposition}
\begin{proof} See, for example, Lemma 1 and Theorem 1 on page 78-79 of \cite{Hs}.\end{proof}

Recall that the Chern Ricci curvature is defined as
$$
\Ric^1(X, Y)=\sum_{i=1}^m \langle \Rm_{X,Y} e_i, \overline{e_i}\rangle=\sqrt{-1}\sum \langle \Rm_{X, Y} \epsilon_i, J\epsilon_i\rangle
$$
where $\{e_i=\frac{1}{\sqrt{2}}(\epsilon_i-J \epsilon_i)\}_{i=1}^m$ is a unitary frame of $T_p^{1, 0}M$.
The easier part of holonomy theorem of Ambrose-Singer asserts that for any $X, Y\in T_p M$, for a metric connection $\nabla$,
$\Rm_{X, Y}\in \mathfrak{g}$ as a skew-symmetric transformation of $T_pM$ (see for example \cite{Ni-21}, Lemma 2.2 for a simple proof of this fact). From this it is easy to see that if $G\subset \mathsf{SU}(m)$, $\operatorname{tr}(\Rm_{X, Y})=0$, hence $\Ric^1(X, Y)=0$.

 The other direction follows from  the harder part of Ambrose-Singer's theorem namely $\{\gamma(\Rm^q)_{X, Y}\}$ generated $\mathfrak{g}$, where for any $q$, $\gamma(\Rm^q)_{X, Y}$ is defined via the idenity $\langle \gamma(\Rm^q)_{X, Y} Z, W\rangle =\langle \Rm_{\gamma^{-1}(X), \gamma^{-1}(Y)}\gamma^{-1}(Z), \gamma^{-1}(W)\rangle$ with $q\in M$, $\gamma$ being the path from $q$ to $p$, ($\gamma$ also denotes the parallel transport along $\gamma$ from $q$ to $p$). Since $\nabla$ is Hermitian, the parallel transport preserves the inner product and the almost complex structure $J$, for $\{e_i\}$ a unitary frame obtained by parallel transport along $\gamma$, we have that
\begin{eqnarray*}
\operatorname{tr}\gamma(\Rm^q)_{X, Y} &=&\sum_{i=1}^m \langle \Rm_{\gamma^{-1}(X),\gamma^{-1}(Y)} \gamma^{-1}(e_i), \overline{\gamma^{-1}(e_i)}\rangle\\
&=&\sum_{i=1}^m \langle \Rm_{\gamma^{-1}(X), \gamma^{-1}(Y)}e_i, \overline{e_i}\rangle\\
&=&\Ric^1_q(\gamma^{-1}(X), \gamma^{-1}(Y))=0\end{eqnarray*}
since  $\Ric^1=0$ at $q$ in $M$.  Thus $\gamma(\Rm^q)_{X,Y}$ has a zero trace.  Since $q$ is arbitrary, $M$ is connected (which is always assumed), and  $\{\gamma(\Rm^q)_{X,Y}\}_{q\in M}$ generates the Lie algebra $\mathfrak{g}$ of $G$, we have that $G\subset \mathsf{SU}(m)$.   In summary we have
\begin{proposition}\label{propo:cy} For a connected Hermitian manifold $(M, g)$ with a Hermitian connection $\nabla$, its restricted holonomy group $G\subset \mathsf{SU}(m)$ if and only if the Chern Ricci $\Ric^1(X, Y)=0$ for any $X, Y\in T_pM$, $\forall p\in M$.
\end{proposition}

\section{Proof of Theorem \ref{thm:CY-torsion}}
We prove the result by establishing a condition under which the center $\mathfrak{z}$ of $\mathfrak{g}$ is nontrivial. The first observation is

\begin{lemma}\label{lemma:31} If $G$ action on $T^{1, 0}M$ is irreducible, then $\mathfrak{z}\subset\mathbb{R}\{\sqrt{-1}\operatorname{id}\}$. In particular it implies that if $\mathfrak{z}\ne 0$, then $\mathfrak{z}=\mathbb{R}\{\sqrt{-1}\operatorname{id}\}$.
\end{lemma}
\begin{proof} The irreducibility implies that for any $a\in \mathfrak{z}$, $a=\lambda \sqrt{-1}\operatorname{id}$ by Proposition \ref{prop:schur}. Note here the almost complex structure $J$ is the same as multiplying by $\sqrt{-1}$ on $T^{1, 0}M$.  It can also be seen by the fact that $a\in\mathfrak{u}(m)$, and on the other hand $a=c\operatorname{id}$ for some complex number $c$ due to the irreducibility and the fact that $[a, \mathfrak{g}]=0$ (this is the Schur's Lemma). This, together with $a+{}^t \bar{a}=0$,  forces $a=\sqrt{-1}\lambda\operatorname{id}$ for some $\lambda\in \mathbb{R}$.
\end{proof}

The second result concludes that the holonomy algebra must contain a nontrivial center if the Chern Ricci $\Ric^1\ne 0$.

\begin{lemma} \label{lemma:32}Let $(M^m, g)$ be a Hermitian manifold with a Hermitian connection $\nabla$. If its Chern Ricci curvature does not vanish, then $\mathfrak{z}\ne \{0\}$, where  $\mathfrak{z}$ is the center of $\mathfrak{g}$.
\end{lemma}
\begin{proof}
Assume otherwise, namely $\mathfrak{z}=\{0\}$, then $\mathfrak{g}=\mathfrak{g}'=[\mathfrak{g}, \mathfrak{g}]$, the derived algebra by Proposition \ref{prop:com-lie-alg}.
Now pick  $a=\lambda \sqrt{-1}\operatorname{id}$ for some real $\lambda\ne 0$. Clearly $[a, \mathfrak{g}]=0$ due to $a$ is in the center of $\mathfrak{u}(m)$. Then we have
\begin{eqnarray*}
0&=&\langle [a, \mathfrak{g}], \mathfrak{g}\rangle\\
&=&\langle a, [\mathfrak{g}, \mathfrak{g}]\rangle\\
&=&\langle a, \mathfrak{g}\rangle.
\end{eqnarray*}
This implies that for any $b\in \mathfrak{g}$, $\langle \sqrt{-1} \operatorname{id}, b\rangle =0$. This then implies that, if  we express $b=x+\sqrt{-1}y$,  $\operatorname{tr}(y)=0$. But $b\in \mathfrak{u}(m)$ also implies that $\operatorname{tr}(x)=0$. Hence we have that $\operatorname{tr}(b)=0$ for any $b\in\mathfrak{g}$. This implies that $G\subset \mathsf{SU}(m)$. We arrive at a contradiction with the assumption that $\Ric^1\ne 0$ by Proposition \ref{propo:cy}.
\end{proof}

Now we  prove Theorem \ref{thm:CY-torsion}. By Lemmata \ref{lemma:31}-- \ref{lemma:32} above we have that $\mathfrak{z}\ne 0$, hence  $J\in \mathfrak{z}\subset \mathfrak{g}$. (Here note that in the identification of $\mathfrak{u}(m)$ with the subalgebra of $\mathfrak{so}(2m)$, multiplying $\sqrt{-1}$ is identified with $J$.)  This implies that
$J(T)=0$ by the assumption that $T$ is parallel. By the definition of the torsion we have that
\begin{eqnarray*}
0&=& (J T)_{X, Y}\\
&=& J(T(X, Y))-T(JX, Y)-T(X, JY)\\
&=&J (\nabla_X Y-\nabla_Y X -[X, Y])-(\nabla_{JX} Y-\nabla_Y JX-[JX, Y])\\
&\quad&-(\nabla_X JY-\nabla_{JY}X-[X, JY])\\
&=&-J[X, Y]-\nabla_{JX} Y+\nabla_{JY}X+[JX, Y]+[X, JY].
\end{eqnarray*}
Here we have used that $[\nabla, J]=0$, namely the connection is Hermitian.
Since we are on a complex manifold, the almost complex structure is integrable implies that Nijenhuis tensor vanishes, which then implies that
$$
JN(X,Y)=J[JX, JY]-J[X, Y]+[X, JY]+[JX, Y]=0.
$$
Combining them together we have that
$$
-\nabla_{JX} Y+\nabla_{JY}X=J[JX, JY].
$$
Applying $J$ again and using that $[\nabla, J]=0$ we have that for any tangent vector field $X, Y$
$$
\nabla_{JX} JY-\nabla_{JY} JX =[JX, JY]
$$
namely $T(JX, JY)=0$. Since $X, Y$ are arbitrary we conclude that $T=0$, hence the manifold must be K\"ahler. Hence we complete the proof of the theorem.

The proof also yields the following result for almost complex manifolds.
\begin{corollary} For an almost complex Hermitian manifold $(M^m, g)$, and any Hermitian connection $\nabla$, assume that  the restricted holonomy group $G$ action on $T^{1, 0}(M)$ is irreducible. Assume  further that $\Ric^{1}\ne  0$ and the torsion is parallel. Then
$$
T(X, Y)=-N(JX, JY)
$$
with $N(\cdot, \cdot)$ being the Nijenhuis tensor of the almost complex structure.
\end{corollary}

When the holonomy action is reducible the proof also gives some result on the torsion restricted to a subbundle which is invariant under the holonomy action. If $\mathcal{V}\subset T^{1, 0}M$ is a subbundle which is irreducible and invariant under the action of the holonomy group $G$. Then $\Rm_{x, y}: \mathcal{V}_p\to \mathcal{V}_p$ for any $x, y$. Let $\Ric^1_{\mathcal{V}}$ denote the Chern Ricci curvature of $\mathcal{V}_p$:
$$
\Ric^1_{\mathcal{V}}(X, Y)=\sum_{i=1}^k \langle \Rm_{X,Y} e_i, \overline{e_i}\rangle
$$
where $\{e_i\}_{i=1}^k$ is a unitary frame of $\mathcal{V}_p$. Let $P^{\mathcal{V}}$ be the projection from $T^{1, 0}_pM$ ($T_pM$) onto $\mathcal{V}_p$ ($\mathcal{V}\oplus \overline{\mathcal{V}}$). By the definition, $P^{\mathcal{V}}$ is also invariant under the parallel transport. We then have the following extension of Theorem \ref{thm:CY-torsion}.

\begin{proposition}\label{prop:lich2}
 For a complex Hermitian manifold $(M^m, g)$, and any Hermitian connection $\nabla$, assume that  the restricted holonomy group $G$ action on $T^{1, 0}(M)$ is reducible with an invariant irreducible subbundle $\mathcal{V}$. Assume  further that $\Ric^{1}_{\mathcal{V}}\ne  0$ and the torsion of $\nabla$ is parallel. Then for any $X, Y\in \mathcal{V}_p\oplus \overline{\mathcal{V}_p}$, $P^{\mathcal{V}}(T(X, Y))=0$.
\end{proposition}
\begin{proof} Without the loss of generality we may assume that $M$ is simply-connected. Now consider $G|_{\mathcal{V}}$ which we denote as $G^{\mathcal{V}}$. The argument in the proof of Proposition \ref{propo:cy} can be adapted to show that  $\Ric^1_{\mathcal{V}}=0$ if any only if $G^\mathcal{V}\subset \mathsf{SU}(\mathcal{V})$.

On the other hand, as in Lemma \ref{lemma:32}, $\Ric^1_{\mathcal{V}}\ne 0$ implies that the Lie algebra $\mathfrak{g}^{\mathcal{V}}$ must have a nontrivial center. This then implies that $J|_{\mathcal{V}}$ belong to $\mathfrak{g}^{\mathcal{V}}$. Hence $J|_{\mathcal{V}} (P^{\mathcal{V}}T)=0$. When $X \in  \mathcal{V}\oplus \overline{\mathcal{V}}$, $J|_{\mathcal{V}} (X)=JX$. Now the above calculation gives the claimed result.
\end{proof}

One may find it interesting to compare the result with the Hopf manifold example in the Appendix of \cite{NZ-CAS2}.

\section{A Lie algebraic structure on $T_\mathbb{C}M$}

In \cite{NZ-TAMS, NZ-CAS2}, the Hermitian manifolds whose Chern connection is Ambrose-Singer was studied. There are several key ingredients in the proof. One of them is that, after splitting off all K\"ahler factors, the manifold $M$ has the property that for any given $p$,  $T^{1, 0}_pM$ (also $T_{\mathbb{C}}M$ admits a complex Lie algebra structure (cf. Lemma 3.1 of \cite{NZ-TAMS}) if the bracket is defined by $[\cdot, \cdot]:=T^c(\cdot, \cdot)$. Here $T^c$ denotes the torsion of the Chern connection. Due to  the assumption  that the torsion is parallel with respect to the Chern connection, it induces a Lie algebraic structure on the tangent bundle $T^{1, 0}M$.

In a recent paper \cite{BPT}, assuming that the Bismut connection has the parallel torsion, and  additionally, the curvature is of the connection is K\"ahler like, the authors exploited the compact Lie algebra structure defined by the torsion of the Bismut connection to obtain a classification result via a result of Samelson and Pittie. Here we prove that for any Hermitian manifold $(M, g)$, if its Bismut torsion is parallel, then the torsion  of the Chern connection endows $T_{\mathbb{C}}M$ with a complex Lie algebra structure.

First we recall the following easy result for the torsion of any Hermitian connection.

\begin{proposition}\label{prop:torsion-Ni}
For any  connection $\nabla$ with $\nabla J =0$, the torsion $T$ satisfies that
\begin{equation}\label{eq:torsion-N}
T_{X, Y}-T_{JX, JY}+J(T_{X, JY}+T_{JX, Y})=N(X, Y), \forall \mbox{ real vectors } X, Y.
\end{equation}
Here $T_{X, Y}$ is $T(X, Y)$.
In particular, on a complex manifold
\begin{equation}\label{eq:torsion-no02}
T_{X+\sqrt{-1}JX, Y+\sqrt{-1}JY}-\sqrt{-1} J(T_{X+\sqrt{-1}JX, Y+\sqrt{-1}JY})=0.
\end{equation}
\end{proposition}
\begin{proof} The first follows by direct calculations:
\begin{eqnarray*}
T_{X, Y}-T_{JX, JY}&=& \nabla_X Y-\nabla_Y X +J(\nabla_{JY} X-\nabla_{JX} Y)+[JX, JY]-[X, Y];\\
T_{X, JY}+T_{JX, Y}&=& J(\nabla_XY-\nabla_YX)-(\nabla_{JY}X-\nabla_{JX} Y)-[X, JY]-[JX, Y]
\end{eqnarray*}
and the definition of $N(X, Y)$. The second follows from the first by expanding the expression and collect terms, noting that $N=0$ for a complex manifold.
\end{proof}

The next result only holds for the Bismut connection.

\begin{theorem}\label{thm:bismut1} Let $(M, g)$ be a Hermitian manifold. Let $\nabla$ be the Bismut connection, which is characterized by that its torsion is  a $3$-form, in the sense that  $\langle T(x, y), z\rangle$ is a $3$ form in the three variables. Assume that $T$ is parallel. Let $T^c$ be the torsion of the Chern connection. Then for any $X, Y\in T^{1, 0}_pM$, $[X, Y]_c=T^c(X, Y)$ defines a Lie bracket, which makes $T^{1, 0}M$ (and $T_{\mathbb{C}}M$) a bundle of complex Lie algebra. Moreover $T_{\mathbb{C}}M=T^{1,0}M\oplus T^{0,1}M$ as a direct sum of two ideals.
\end{theorem}
\begin{proof} It is well known that the torsion $T^c$ of the Chern connection is a map of   $\wedge^2T^{1, 0}M$ into $T^{1, 0}M$. Namely $T^c(X, \overline{Y})=0=T^c(\overline{X}, Y)$ for any $X, Y$ vectors of $(1, 0)$ type. It is also a map of $\wedge^2T^{0,1}M$ into $T^{0,1}M$. The torsion of $T$ and $T^c$ are related by the following formulae (cf. \cite{LS}):
\begin{equation}\label{eq:two-ts}
\langle T_{X, Y}, \overline{Z}\rangle =-\langle T^c_{X, Y}, \overline{Z}\rangle;\,  \langle T_{\overline{X}, Y}, \overline{Z}\rangle =\langle \overline{T^c_{X, Z}}, Y\rangle; \,  \langle T_{\overline{X}, Y}, Z\rangle =-\langle T^c_{Y, Z}, \overline{X}\rangle.
\end{equation}
One can derive the rest of $T$ by the total screw symmetry. The above proposition asserts that
$$
\langle T_{X, Y}, Z\rangle=0, \forall X, Y, Z \in T^{1,0}M \mbox{ or } X, Y, Z \in T^{0,1}M.
$$
Namely $T(\wedge^2 T^{1, 0}M )\subset T^{1, 0}M$ and $T (\wedge^2 T^{0, 1}M) \subset T^{0, 1}M$. On the other hand, by the generalized first Bianchi identity:
$$
  {\mathfrak S}\{ (\nabla_xT)(y,z) - \Rm_{xy}z - T(x, T(y,z)) \} = 0
$$
we have in particular, since $\nabla T=0$,
 \begin{equation}\label{eq:B1-2}
  {\mathfrak S}\{  \Rm_{X, \bar{Y}}Z + T(X, T(\overline{Y},Z)) \} = 0.
 \end{equation}
 Here ${\mathfrak S}$ is the cyclic permutation operator of three variables on the expression. On the other hand, the curvature of any Hermitian connection satisfies that $\Rm_{X, \bar{Y}} (T^{1, 0}M) \subset T^{1, 0}M$ and
  $\Rm_{X, \bar{Y}} (T^{0, 1}M) \subset T^{0, 1}M$. Moreover the parallelness of the Bismut torsion (cf. \cite{ZZ22}) implies that
  $$
  \Rm_{Z, X}=0, \forall Z, X \in T^{1,0}M\footnote{The proof,  effectively shows that $ {\mathfrak S}\{\Rm_{X, Y}Z\}=0$ implies that $\Rm_{Z, X}=0$, namely the existence of the Lie algebra structure is equivalent to $\Rm_{Z, X} =0$, since $\Rm_{Z, X}\overline{Y} =0$ is enough to imply $\Rm_{Z, X}=0$.}.
  $$
 Thus by checking the projection of (\ref{eq:B1-2}) into $T^{0,1}M$ we have that  for any $X, Y, Z\in T^{1, 0}M$
 $$
 P_2( T_{X, T_{\overline{Y}, Z}}+T_{\overline{Y}, T_{Z, X}}+T_{Z, T_{X, \overline{Y}}})=0.
 $$
 Here $P_2$ ($P_1$) denotes the projection into $T^{0,1}M$ ($T^{1, 0}M$). Using the proposition and (\ref{eq:two-ts}) we have that (here we use the Einstein's convention of summing the repeated indices)
 \begin{eqnarray*}
T_{\overline{Y}, Z}&=& -\langle T^c_{Z, e_k}, \overline{Y}\rangle \overline{e_k}+I;\\
T_{Z, X}&=& - T^c_{Z, X};\\
T_{X, \overline{Y}}&=& \langle T^c_{X, e_k},  \overline{Y}\rangle \overline{e_k}+II.
 \end{eqnarray*}
 Here $I$ and $II$ denote the components of the corresponding left hand side  in $T^{1, 0}M$ (namely $P_1(T_{\overline{Y}, Z})$ and $P_1(T_{X, \overline{Y}})$). Since $P_2(T_{X, I})=P_2(T_{Z, II})=0$ we have that
 \begin{eqnarray*}
 0&=&P_2(-T_{X, \overline{e_k}}\langle T^c_{Z, e_k}, \overline{Y}\rangle- T_{\overline{Y}, e_k}\langle T^c_{Z, X}, \overline{e_k}\rangle+T_{Z, \overline{e_k}}\langle T^c_{X, e_k},  \overline{Y}\rangle)\\
 &=&\left(-\langle T^c_{Z, e_k}, \overline{Y}\rangle \langle T^c_{X, e_\ell}, \overline{e_k}\rangle+
\langle T^c_{Z, X}, \overline{e_k}\rangle\langle T^c_{e_k, e_\ell}, \overline{Y}\rangle+ \langle T^c_{X, e_k},  \overline{Y}\rangle\langle T^c_{Z, e_\ell}, \overline{e_k}\rangle \right) \overline{e_\ell}\\
&=&-\langle T^c_{Z, T^c_{X, e_\ell}}+T^c_{e_\ell, T^c_{Z, X}}+T^c_{X, T^c_{e_\ell, Z}}, \overline{Y}\rangle  \overline{e_\ell}.
\end{eqnarray*}
This implies that
$$
 T^c_{Z, T^c_{X, e_\ell}}+T^c_{e_\ell, T^c_{Z, X}}+T^c_{X, T^c_{e_\ell, Z}}=0, \forall \ell.
$$
This checks that the Lie bracket defined by the torsion of the Chern connection satisfies the Jacobi identity.
By (\ref{eq:two-ts}) we also have that $\nabla T^c=0$. Hence $T^{1, 0}M$ is a bundle of Lie algebra.

Since for the Chern connection $T^c(X, \overline{Y})=0, \forall X, Y \in T^{1, 0}M$, the Lie bracket structure on $T^{1,0}M$ and $T^{0, 1}M$ defined by $T^c$ extends to one on $T_{\mathbb{C}}M$.
\end{proof}

Since we do not know if $\nabla^c T^c=0$, the result may  not be derived from the generalized first Bianchi identity for the Chern curvature and Chern torsion as in \cite{NZ-TAMS}. However one can derive the following consequence via the generalized 1st Bianchi identity.

\begin{corollary}
Let $(M, g)$ be a Hermitian manifold. Let $\nabla$ be the Bismut connection. Assume that $T$ is parallel.
 Then for any $X, Y\in T^{1, 0}_pM$, $[X, Y]=T(X, Y)$ defines a Lie bracket, which makes $T^{1, 0}M$ a bundle of complex Lie algebra.
\end{corollary}

\section{Proof of Theorem \ref{thm:main1}}

Since for any given smooth function $\varphi$,  adding $\epsilon \sqrt{-1}\partial \bar{\partial} \varphi$ with a small $\epsilon$ to a K\"ahler metric with $\mathsf{U}(m)$ holonomy (with respect to the Levi-Civita connection) does not change the holonomy being $\mathsf{U}(m)$, one may say that  generic K\"ahler metrics (manifolds) have $\mathsf{U}(m)$ holonomy. The proof of Theorem \ref{thm:main1} follows  the line of argument in \cite{Beau} together with tracing the related estimates for the positive (quasi-positive) case. A  well-known result of Kodaira (Theorem 8.3 of \cite{MK}), which asserts that {\it a K\"ahler manifold  $M$ is projective if the Hodge number $h^{2, 0}=\dim (\mathcal{H}^{2, 0}(M))$, the dimension of the harmonic $(2, 0)$-forms, vanishes}, reduces the issue to the vanishing of $h^{2, 0}(M)$.

\begin{proposition}\label{prop:1} Assume that $(M, g)$ is a compact K\"ahler manifold. Let $\lambda_i(x)$ be the eigenvalue of its Ricci form. Assume that
\begin{equation}\label{eq:p-sum}
\min_{1\le i_1<i_2<\cdots<i_p\le m} \left(\lambda_{i_1}+\cdots +\lambda_{i_p}\right)\ge 0.
\end{equation}
Then any harmonic $(p, 0)$-form is parallel.
\end{proposition}
\begin{proof} Let $\phi$ be a harmonic $(p, 0)$-form. Since it is automatically holomorphic, we can apply the well-known Bochner formula (cf. Lemma 3.1 of \cite{NST} and \cite{MK}) which asserts that
\begin{equation}\label{eq:key}
 -\square_{\bar{\partial}} |\phi|^2\ge  |\nabla \phi|^2 + \min_{1\le i_1<i_2<\cdots<i_p\le m} \left(\lambda_{i_1}+\cdots +\lambda_{i_p}\right) |\phi|^2
\end{equation}
where $|\phi|^2$ is the square of the norm of $\phi$ and $\square_{\bar{\partial}}=\bar{\partial} \bar{\partial}^*+\bar{\partial}^*\bar{\partial}$ is the Hodge-Laplacian. Integrating (\ref{eq:key}) on $M$ we have that $\nabla \phi\equiv 0$, namely $\phi$ is parallel on $M$.
\end{proof}
In particular $\phi$ is invariant under the parallel transport along any closed piece-wisely smooth path, hence the holonomy group action. Thus a nonzero parallel $\phi$ implies the existence of a one dimensional invariant subspace of $\wedge^p T'_{x_0}M$ under  the action of the relative holonomy group $H^0_{x_0}$.  Now apply the  result  below from the representation theory (cf. Corollary 5.5.3  of \cite{GW} and exercise 5.6.1 of \cite{BD}) we can conclude that any parallel $(p, 0)$ form $\phi\equiv0$, for $1\le p< m$.

\begin{proposition}\label{prop:rep} The action of $\mathsf{U}(m)$ ($\mathsf{SU}(m)$) on $\wedge^p\mathbb{C}^m$ is irreducible. Same result holds on $\otimes^p \mathbb{C}^m$.
\end{proposition}
\begin{proof} We include a short proof of this  result for the convenience of people without any background of Lie theory.
Consider  $(\rho, V )$ a finite dimensional irreducible unitary representation of $\mathsf{U}(m)$ on $ V = \wedge^p \mathbb{C}^m$ with $1\le p <m$ and $\rho (g) = \wedge^p g$. Let T denote the diagonal elements of $\mathsf{U}(m)$ (a maximal torus of $\mathsf{U}(m)$) and let $e_1,\cdots, e_m$ be the standard basis of $\mathbb{C}^m$. If $I_p\subset \{1, \cdots, m\}$  is a subset with $p$ elements, $1 \le i_1 < \cdots < i_p\le m$ be the elements then set $e_{I_p} = e_{i_1}\wedge\cdots \wedge e_{i_p}$. If $z\in T $ has diagonal entries $z_1, \cdots, z_m$ then $\rho(z)(e_{I_p}) = \Pi_{j=1}^p z_{i_j} e_{I_p}$. Set $\chi_{I_p}(z) = \Pi_{j=1}^p z_{i_j}$. Since the $e_{I_p}$ form a basis of $V$ we see that if $v \in V$ expressed as $v=\sum_{I_p} v_{I_p}\cdot  e_{I_p}$ then
\begin{equation}\label{eq:51}\rho(z)\, v =\sum_{I_p} \chi_{I_p}(z)\cdot v_{I_p}\cdot e_{I_{p}}.\end{equation} Integrating this over $T$ with respect to the normalized Haar measure $d\mu$, direct calculation shows that (as in the Fourier expansion)
\begin{equation}\label{eq:52} v_I\cdot e_I =\int_{T}\chi^{-1}_{I}(z)\rho(z)\,  v d\mu(z).\end{equation} Thus if $W$, with $0 \subsetneq W \subset V$,  is a irreducible invariant subspace then at least one of the $e_I \in  W$. Indeed since there exists some $v\ne 0$, $v\in W$, $\rho(z)v \in W$ for all $z$. The equation (\ref{eq:52}) implies that $v_I\cdot e_I$ all in $W$. Since they can not all vanish, there must be one $ v_I\cdot e_I\ne 0$ belongs to $W$.  Let $S$ be the subgroup of permutation matrices in $\mathsf{U}(m)$ that is, if $s\in S_m$ then $s(e_i) = e_{s(i)}$ is the corresponding element of $S$. We note that $\rho(s)\, e_I = e_{s(I)} = e_{i_{s(1)}} \wedge \cdots \wedge e_{i_{s(p)}}$. Thus if $W$ is a nonzero invariant subspace of V then $e_{I_p} \in  W$ for all $I_p\subset \{1, \cdots, m\}$ with $p$ elements. Thus $W = V$. This proves that induced action of $\mathsf{U}(m)$ on $\wedge^p \mathbb{C}^m$ is irreducible. \end{proof}

Combining Proposition \ref{prop:1} and Proposition \ref{prop:rep} we have a proof for the part (i) of Theorem \ref{thm:main1}.

\begin{corollary}\label{coro:1} Assume that $(M, g)$ satisfies (\ref{eq:p-sum}) and with holonomy $\mathsf{U}(m)$. Then $h^{p', 0}=0$ for all $p\le p'\le m$.
\end{corollary}

Now we extend the result to other curvatures. We first deal with the mixed curvature case. We include a simple consequence of positive holomorphic sectional curvature, even though it is not used directly.

\begin{lemma}\label{lemma:51} If at a point $x_0\in M$, the holomorphic sectional curvature $H(X)\ge \kappa |X|^4$, then there exists a unitary frame $\{e_i\}_{i=1}^m$ such that
\begin{equation}\label{eq:53}
\Ric(e_i, \overline{e}_i)\ge \frac{\kappa}{2}(m+1).
\end{equation}
\end{lemma}
\begin{proof} Let $e_1$ be the unit direction where $H$ attains its minimum $\kappa$. The second variational argument of \cite{N} (cf. Corollary 2.1 of \cite{NZ}) implies that
$$
\Rm(e_1, \overline{e}_1, X, \overline{X})\ge \frac{\kappa}{2}, \forall X \perp e_1, |X|=1.
$$
Now we pick $e_2\in (e_1)^{\perp}$ such that $H(e_2)=\inf_{X\in (e_1)^\perp, |X|=1} H(X)$. Clearly $H(e_2)\ge \kappa$. By the same second variation argument we have that
$$
\Rm(e_2, \overline{e}_2, X, \overline{X})\ge \frac{\kappa}{2}, \forall X \perp \{e_1, e_2\}, |X|=1.
$$
Note that we do have $\Rm(e_2, \overline{e}_2, e_1, \overline{e}_1)\ge \frac{\kappa}{2}$ by the first step.
Repeat the above procedure, we find the frame $\{e_i\}$. From its construction we have that
$$
\Ric(e_i, \overline{e}_i)\ge \kappa +(m-1)\frac{\kappa}{2}=\frac{\kappa}{2}(m+1).
$$
This proves the claim.
\end{proof}

However, it does not imply that $\Ric\ge \frac{\kappa}{2}(m+1)$. Now we prove part (ii) of the theorem.

Recall the Bochner formula for a holomorphic $(p, 0)$ form which is viewed as a holomorphic section of $\wedge^p (T^{1, 0}M)^*$.

\begin{lemma}\label{lemma:52} Let $s$ be a global holomorphic $p$-form on $M^m$ which locally is expressed as  $s=\frac{1}{p!}\sum_{I_p} f_{I_p}\varphi_{i_1}\wedge \cdots \wedge \varphi_{i_p}$ ,  where $I_p=(i_1, \cdots, i_p)$ and $\{ \varphi_1, \ldots , \varphi_m\}$ is a local unitary coframe. Then
  $$\partial \overline{\partial } \,|s|^2 = \langle \nabla s , \overline{\nabla s} \rangle  - \widetilde{R}(s,
  \overline{s}, \cdot , \cdot ) $$
where $\widetilde{R}$ stands for the curvature of the Hermitian bundle $\bigwedge^p\Omega$, where $\Omega=(T'M)^*$ is
the holomorphic cotangent bundle of $M$. The metric on $\bigwedge^p\Omega$ is derived from the metric of $M^m$. Then for any unitary coframe $\{\varphi_i\}$,
\begin{equation}\label{eq:54}
\langle \sqrt{-1}\partial\bar{\partial} |s|^2, \frac{1}{\sqrt{-1}}v\wedge \bar{v}\rangle =\langle \nabla_v s,
\overline{\nabla_vs}\rangle +\frac{1}{p!}\sum_{I_p} \sum_{k=1}^p  \sum_{l=1}^m  R_{v\bar{v}i_k \bar{l}}\,f_{I_p} \overline{f_{i_1\cdots(l)_k\cdots i_p}}.
\end{equation}
 Also, given any $x_0$, there exists a unitary coframe $\{\varphi_i\}$ at $x_0$, such that
\begin{equation}\label{eq:55}
\langle \sqrt{-1}\partial\bar{\partial} |s|^2, \frac{1}{\sqrt{-1}}v\wedge \bar{v}\rangle =\langle \nabla_v s,
\overline{\nabla_v s}\rangle +\frac{1}{p!}\sum_{I_p} \sum_{k=1}^p R_{v\bar{v}i_k \bar{i}_k}|f_{I_p}|^2.
\end{equation}
\end{lemma}
\begin{proof} One can refer to Lemma 2.1 of \cite{NZ-2}.  Here one simply diagonalize the Hermitian form $f_{I_p} \overline{f_{i_1\cdots(l)_k\cdots i_p}}$ in terms of $(i_k, \bar{l})$.
\end{proof}

We also need the following formula from \cite{Zhang-Zhang} (Lemma 3.2).

\begin{lemma}[Zhang-Zhang]\label{L1}
		Let $(M, g)$ be a Hermitian manifold of complex dimension $m$, $\eta$ be a real $(1,1)$-form  and $s$ be a $(p,0)$-form on $M$. We can define a real semi-positive $(1,1)$-form $\beta$ associated to $s$ by
		\begin{equation}\label{eq:56}
			\widetilde{\beta}=\Lambda{}^{p-1}\left((\sqrt{-1})^{p^2}\frac{s\wedge\bar{s}}{p!}\right),
		\end{equation}
		where $\Lambda $ is the contraction operator dual to the Lefschetz operator $L$, namely  wedge product by the K\"ahler form $ \omega=\sqrt{-1}g_{i\bar{j}}dz^i\wedge d\bar{z}^j$. Then
		\begin{equation}\label{eq:57}
			\operatorname{tr}_\omega \widetilde{\widetilde{\beta}}=\Lambda \widetilde{\beta}=\vert s\vert^2_g
		\end{equation}
		and
		\begin{equation}\label{eq:58}
			\begin{split}
				(\sqrt{-1})^{p^2}\eta\wedge s\wedge\bar{s}\wedge \frac{\omega^{m-p-1}}{(m-p-1)!}=\left[\tr_\omega\eta\cdot\vert s\vert_g^2-p\langle \eta,\widetilde{\beta}\rangle_g\right]\frac{\omega^m}{m!}.
			\end{split}
		\end{equation}
\end{lemma}		

Now we may adapt the argument of \cite{NZ, NZ-2, Zhang-Zhang, Tang} to prove that any harmonic (hence holomorphic) $(2, 0)$-form $\phi$ is invariant under the holonomy action, hence vanish by Proposition \ref{prop:rep}.

As in \cite{NZ, NZ-2} applying the above to $s=\wedge{}^k\phi$ which is a holomorphic $(2k, 0)$ form, where $k$ is the largest integer satisfying that $\wedge{}^k \phi\ne 0$ and $\wedge{}^{k+1}\phi=0$.  The key is that after choosing a suitable unitary frame $\phi=\sum_{i=1}^k f_i dz^{2i-1}\wedge dz^{2i}$. This makes $s$ a simple $(2k, 0)$-form. The rest of the  argument runs closely as that of \cite{Tang}.

First by Lemma \ref{lemma:52}, for the unitary frame $\{e_i\}_{i=1}^m=\{\frac{\partial}{\partial z^i}\}_{i=1}^m$ chosen to put $\phi$ into a normal form and $s$ into a simple $(2k, 0)$ form, we have that
\begin{eqnarray}
\Delta |s|^2&=&|\nabla s|^2+|s|^2 \sum_{i=1}^{2k} \Ric_{i\bar{i}};\label{eq:59}\\
\langle \sqrt{-1} \partial\bar{\partial} |s|^2, \widetilde{\beta}\rangle &=&\frac{|s|^2}{2k}\sum_{i=1}^{2k}|\nabla_i s|^2+\frac{|s|^4}{2k}\sum_{i, j=1}^{2k} \Rm_{i\bar{i}j\bar{j}}. \label{eq:510}
\end{eqnarray}
Note that $\widetilde{\beta}=\frac{|s|^2}{2k}\left(dz^1\wedge d\bar{z}^1+\cdots+dz^{2k}\wedge d\bar{z}^{2k}\right)$.
Since $s$ is holomorphic, applying Lemma \ref{L1} to the case $\eta=\sqrt{-1}\partial\bar{\partial}|s|^2$ implies that
\begin{equation}\label{eq:511}
\int_M \left(\Delta |s|^2\right) |s|^2 \omega^m=2k\int_M \langle \sqrt{-1}\partial \bar{\partial} |s|^2, \widetilde{\beta}\rangle \omega^m.
\end{equation}
Note that by integrating (\ref{eq:510}) on $M$, we have that the right hand side above is nonnegative if the $2k$-scalar curvature $S_{2k}\ge 0$.

By the assumption that the mixed curvature is nonnegative, a Berger's averaging trick (cf. Lemma 1.1 of \cite{NZ-2}, Appendix of \cite{Ni-1807}),  implies that
$$
\frac{\alpha}{2k}\sum_{i=1}^{2k}\Ric_{i\bar{i}}+\frac{\beta}{k(2k+1)}\sum_{i, j=1}^{2k}\Rm_{i\bar{i}j\bar{j}}\ge 0.
$$
Together with the above consequence of Lemma \ref{L1} and
integrating (\ref{eq:59}), after multiplying $|s|^2$ on the both sides,  and (\ref{eq:510}) over $M$ and linearly combining them as in \cite{Tang},  we have that
\begin{eqnarray*}
\left( \frac{\alpha}{2k}+\frac{\beta}{k(2k+1)}\right)\int_M \Delta |s|^2 \cdot|s|^2 \omega^m &\ge& \int_M |s|^2 \left(\frac{\alpha}{2k}|\nabla s|^2+\frac{\beta}{k(2k+1)}\sum_{i=1}^{2k}|\nabla_i s|^2\right).
\end{eqnarray*}
Since the left hand side is non-positive, by the assumption on $\alpha$ and $\beta$, we have that $|\nabla s|=0$. Now Proposition \ref{prop:rep} implies that $s$  vanishes. This leads to a contradiction that $\phi$ is a nontrivial holomorphic $(2, 0)$-form. Hence we have $h^{2, 0}(M)=0$.

The above argument for $\mathcal{C}_{\alpha, \beta}\ge 0$ is basically that of \cite{Tang}, which was motivated by \cite{Zhang-Zhang}\footnote{There are also two more recent papers \cite{CL, Matsumura} studying the fundamental groups of manifolds with nonnegative holomorphic  curvature.}. The $S_2\ge 0$ case requires a bit new ingredients however. First by  (\ref{eq:510}) and (\ref{eq:511}), since $S_2\ge 0$ implies that $S_{2k}\ge 0$, namely $\sum_{i, j=1}^{2k} \Rm_{i\bar{i}j\bar{j}}\ge 0$ which implies that the right hand side of (\ref{eq:510}) is nonnegative, we conclude that if $s$ is nonzero,
$\nabla_i s=0$ for $1\le i\le 2k$, $|s|^2$ is a constant (from the integration by part of the left hand side of (\ref{eq:511})) and $\sum_{i, j=1}^{2k}\Rm_{i\bar{i}j\bar{j}}$ vanishes.  Hence at the $\Sigma^{2k}=\{ \frac{\partial}{\partial z^{i}}\}_{i=1}^{2k}$ , $S_{2k}( \Sigma')$ attains its minimum ($0$) among all $2k$-dimensional subspaces. Now by the second variation argument (cf. Proposition 4.2 of \cite{N}) we have that
$$
\sum_{i=1}^{2k}\Rm_{i\bar{i} \ell \bar{\ell}}\ge 0, \forall \ell\ge 2k+1.
$$
  Now apply (\ref{eq:54}) to $v=\ell$  we can conclude that $|\nabla_{\ell} s|^2=0$ for $\ell\ge 2k+1$. Thus $s$ is parallel. We then have proved the  claimed result in part (iii).

\begin{remark}
1. In \cite{CDP}  it was shown that if the holonomy group of a projective manifold is $\mathsf{U}(m)$, the manifold is rationally connected. Hence the manifold in  Theorem \ref{thm:main1} is rationally connected.

2. If $M$ is simply-connected,  the reducibility of the holonomy action on $T^{1, 0}M$ will split the manifold as products of K\"ahler manifolds. Since the curvature assumptions of Theorem \ref{thm:main1} pass down to the factors, by Berger's classification of the Riemannian holonomy applied to the K\"ahler setting,  one still has results for the factors whose holonomy are either $\mathsf{U}(k)$ or $\mathsf{SU}(k)$, with $k$ being the complex dimension of the corresponding factor.
\end{remark}

\section{Appendix}

Y.-Y. Niu found that one can derive Lemma \ref{L1} of \cite{Zhang-Zhang} using standard formulae involving the the operators $L$, $\Lambda$ and their relation with the Hodge star $*$ operator in \cite{Demailly}. The computation recorded below is due to her.

Using the notation above we have that
\begin{eqnarray*}
\langle \alpha \wedge s, \omega\wedge s\rangle&=& \langle  \alpha \wedge s, L(s)\rangle\\
&=&\langle \Lambda (\alpha \wedge s), s\rangle\\
&=& \operatorname{tr}_\omega (\alpha)|s|^2-p\langle\alpha, \beta\rangle.
\end{eqnarray*}
On the other hand, direct calculations show that
\begin{eqnarray*}
\langle \alpha \wedge s, \omega\wedge s\rangle \frac{\omega^m}{m!}&=& \alpha \wedge s\wedge \overline{*L(s)}\\
&=&\alpha \wedge s\wedge \overline{\Lambda *(s)}\\
&=&\frac{(\sqrt{-1})^{p^2}}{(m-p)!}\alpha \wedge s \wedge\overline{\Lambda\cdot L^{m-p}(s)}\\
&=& (\sqrt{-1})^{p^2}\eta\wedge s\wedge\bar{s}\wedge \frac{\omega^{m-p-1}}{(m-p-1)!}.
\end{eqnarray*}
Here  at the second line we used $L*=*\Lambda$, at the third line we used calculation of $*s$,  at the last line we used that $\Lambda s=0$ and the formula $[L^\ell,\Lambda](\alpha)=\ell (k-m+\ell-1) L^{\ell-1}(\alpha)$ for any $k$-form $\alpha$.
This proves Lemma \ref{L1}.

\section*{Acknowledgments} {} The author would like to thank S. Sam and N. Wallach for communicating the proof of Proposition \ref{prop:rep} during the peak of the pandemic in 2020,  Y.-Y. Niu for bringing the preprint \cite{CL} to his attention and a careful reading of the first draft. He also thanks R. Lafuente and University of Queensland for the hospitality.


\begin{thebibliography}{A}


\bibitem{AS} W. Ambrose and I.M. Singer,  \textit{ On homogeneous Riemannian manifolds.} Duke Math. J. \textbf{25} (1958), 647--669.


\bibitem{AndV} A. Andrada, R. Villacampa,  \textit{ Bismut connection on Vaisman manifolds.} Math. Zeit. \textbf{302} (2022), 1091-1126.

\bibitem{AOUV} D. Angella, A.Otal, L. Ugarte and R. Villacampa, \textit{
 On Gauduchon connections with K\"ahlerh-like curvature.} Comm. Anal. Geom. \textbf{30} (2022), no. 5, 961--1006.



\bibitem{Berger} M. Berger, \textit{ Sur les groupes d'holonomie homog\`ene des vari\'et\'es \`a connexion affine et des vari\'et\'es riemanniennes.}  Bull. Soc. Math. France \textbf{83} (1955),
    279--330.

\bibitem{Besse} A. L. Besse, \textit{ Einstein Manifolds.} Reprint of the 1987 edition. Classics in Mathematics. Springer-Verlag, Berlin, 2008. xii+516 pp. ISBN: 978-3-540-74120-6.

\bibitem{Beau} A. Beauville, \textit{ Vari\'et\'es k\"ahl\'eriennes compactes avec $c_1=0$.}  Geometry of K3 surfaces: moduli and periods (Palaiseau, 1981/1982). Ast\'erisque No. \textbf{126}(1985), 181--192.

    \bibitem{Boothby} W. Boothby, \textit{ Hermitian manifolds with zero
curvature.} Michigan Math. J., \textbf{5} (1958), no. 2, 229--233.

\bibitem{BD} T.   Br\"ocker and T.  Dieck, \textit{ Representations of compact Lie groups.}  Springer-Verlag, New York, 1985. x+313 pp.

 \bibitem{BPT} G. Barbaro, F. Pediconi and N. Tardini, \textit{ Puriclosed manifolds with parallel Bismut trosion.} ArXiv: 2406.070339.


\bibitem{CDP} F.  Campana, J.-P.  Demailly and T. Peternell, \textit{ Rationally connected manifolds and semipositivity of the Ricci curvature.} Recent advances in algebraic geometry, 71--91, London Math. Soc. Lecture Note Ser., 417, Cambridge Univ. Press, Cambridge, 2015.


\bibitem{CLT}    J. Chu, M. Lee and L.F. Tam,\textit{ K\"ahler manifolds and mixed curvature.}
Trans. Amer. Math. Soc. \textbf{375} (2022), no. 11, 7925--7944.

\bibitem{CL}    J. Chu, M. Lee and J. Zhu,\textit{ On  K\"ahler manifolds with nonnegative mixed curvature.}
ArXiv:2408.14043.

\bibitem{Demailly} J.-P.Demailly, \textit{Analytic Methods in Algebraic Geometry.} High educational press, Beijing and International Press, Boston, 2010.

\bibitem{GHR} S.J. Gates, C.M. Hull and M. Roc\u{e}k, \textit{ Twisted multiplets and new supersymmetric nonlinear sigma models.} Nuc. Phys. B \textbf{248} (1984), 157-186.


\bibitem{GW} R.   Goodman and N. Wallach, \textit{ Symmetry, representations, and invariants.}  Springer, Dordrecht, 2009. xx+716 pp.

\bibitem{Helgason} S. Helgason, \textit{ Differential geometry, Lie groups, and symmetric spaces.} American Mathematical Society, Providence,
RI, 2001.

\bibitem{Hs} W. Hsiang, \textit{ Lectures on Lie Groups.} World Scientific, 2000.

\bibitem{KN} S. Kobayashi, K. Nomizu, \textit{ Foundations of Differential Geometry.} Vol. 1, Interscience Publisher, New York, 1963.

\bibitem{Kostant-TAMS} B. Kostant, \textit{ Holonomy and the Lie algebra of infinitesimal motions of a Riemannian manifold.}  Trans. Amer. Math. Soc. \textbf{80} (1955), 528--542.

    \bibitem{Kostant-Nagoya}  B.   Kostant, \textit{ A characterization of invariant affine connections.} Nagoya Math. J. \textbf{16} (1960), 35--50.

\bibitem{LS}  R.  Lafuente and J. Stanfield,  \textit{H ermitian manifolds with flat Gauduchon connections.} Scuola Normale Superiore di Pisa. Annali. Classe di Scienze, \textbf{11}.

    \bibitem{MK}  J.   Morrow and K.  Kodaira, \textit{ Complex manifolds.} Holt. Rinehart and Winston, New
     York-Montreal-London, 1971.


 \bibitem{Matsumura} S.-I Matsumura, \textit{Fundamental groups of compact K\"ahler manifolds with semi-positive holomorphic sectional curvature.}   Preprint.

\bibitem{Ni-1807}L. Ni, \textit{ Liouville  theorems and a Schwarz Lemma for holomorphic mappings between K\"ahler manifolds.} Comm. Pure Appl. Math., \textbf{74} (2021), 1100--1126.

\bibitem{Ni-21} L. Ni, \textit{ An alternatie induction argument in Simons' proof of holonomy theorem.} Analysis and Partial Differential Euqations on Manifolds, Fractals and Graphs (Nankai, 2019, A. Grigoryan, Y Sun, Eds.) Advances in Analysis and Geometry, \textbf{3} (2021), 443--458.

\bibitem{N} L. Ni, \textit{ The fundamental group, rational connectedness and the positivity of K\"ahler manifolds.} J. reine angew. Math.(Crelle), \textbf{774}(2021), 267--299. DOI 10.1515/crelle-2020-0040.



\bibitem{NST} L. Ni, Y.  Shi and L.-F.  Tam, \textit{ Poisson equation, Poincar\'e-Lelong equation and curvature decay on complete K\"ahler manifolds.} J. Differential Geom. \textbf{57} (2001), no. 2, 339--388.


\bibitem{NZ} L. Ni and F. Zheng, \textit{ Comparison and vanishing  theorems for K\"ahler manifolds.} Calc. Var. Partial Differential Equations, \textbf{57}(2018), no. 6, Art. 151, 31 pp.

\bibitem{NZ-2} L. Ni and F. Zheng, \textit{ Positivity and Kodaira embedding theorem.} Geom. Topol. \textbf{26} (2022), no. 6, 2491--2506.

\bibitem{NZ-TAMS}L. Ni and F. Zheng, \textit{ On Hermitian manifolds whose Chern connection is Ambrose-Singer.}  Trans. Amer. Math. Soc. \textbf{376} (2023), no.9, 6681--6707.

\bibitem{NZ-CAS2} L. Ni and F. Zheng,  \textit{ A classifocation of locally Chern homogenous Hermitian manifolds.} ArXiv preprint.



\bibitem {Stro} A. Strominger, \textit{ Superstrings with torsion.} Nuclear Phys. B. \textbf{274} (1986), 253--284.

\bibitem{Tang} K. Tang, \textit{ Quasi-positive curvature and vanishing theorems.} ArXiv:2405.03895



\bibitem{Vez} L. Vezzoni, \textit{ A note on canonical Ricci forms on $2$-step nilmanifolds.(English summary)} Proc. Amer. Math. Soc. \textbf{141}(2013), no.1, 325--333.


 \bibitem{Wilking} B. Wilkng, \textit{ A Lie algebraic approach to Ricci flow invariant curvature conditions and Harnack inequalities.} J. Reine Angew. Math. \textbf{679} (2013), 223--247.



\bibitem{Wu}  H. Wu and W.-H. Chen, \textit{ Selected Topics on Riemannian Geometry (in chinese)}. Peking Univ. Press, Beijing, 1993.

\bibitem{Yang}
		X. Yang,\textit{ RC-positivity, rational connectedness and Yau's conjecture.}
		 Cambridge J.  Math. \textbf{ 6} (2018), No. 2, 183--212.
		



 \bibitem{Yau-op} S.-T. Yau,\textit{ Open problems in geometry.} Differential Geometry: partial differential equations on manifolds, Proc. Sympos. Pure. Math., vol. \textbf{54}, Amer. Math. Soc., Providence, RI, 1993, pp. 1--28.

     \bibitem{Zhang-Zhang} S. Zhang and X. Zhang, \textit{ On the structure of compact K\"ahler manifolds with nonnegative holomorphic sectional curvature.} ArXiv:2311.18779.

 \bibitem{Zhao-Zheng} Q.   Zhao and F. Zheng, \textit{
Strominger connection and pluriclosed metrics.} (English summary)J. Reine Angew. Math. \textbf{796} (2023), 245--267.

\bibitem{ZZ22} Q. Zhao, F. Zheng, \textit{ On Hermitian manifolds with Bismut-Strominger parallel torsion.} ArXiv: 2208.03071.

\end{thebibliography}
\end{document}